\documentclass{IEEEtran}

\usepackage{amsthm}
\usepackage{epsfig} % for postscript graphics files
\usepackage{epstopdf}
\usepackage{amsmath}
\usepackage{comment}
\usepackage{url}
\usepackage{algorithm,float}
\usepackage{algpseudocode}

\usepackage{amsfonts}
\usepackage{bm}
\usepackage{color}
\usepackage{cite}
\usepackage{stfloats}
\usepackage{chemarrow}
\usepackage{extarrows}

\usepackage{url}
\usepackage[switch]{lineno}
\usepackage{multirow}

\usepackage{setspace}

\newboolean{showcomments}
\setboolean{showcomments}{true}
\newcommand{\wang}[1]{\ifthenelse{\boolean{showcomments}}
	{ \textcolor[rgb]{1,0,1}{(ZW:  #1)}}{}}
\newcommand{\Xie}[1]{\ifthenelse{\boolean{showcomments}}
	{ \textcolor{blue}{(YX:  #1)}}{}}
\newcommand{\jpang}[1]{\ifthenelse{\boolean{showcomments}}
	{\textcolor[rgb]{0,0,0}{#1}}{}}

\theoremstyle{definition}
\newtheorem{theorem}{Theorem}
\newtheorem{lemma}{Lemma}

\theoremstyle{definition}

\newtheorem{remark}{Remark}

\newtheorem{assumption}{\textit{Assumption}}

\allowdisplaybreaks

\begin{document}
	
	\title{Distributed Online Generalized Nash Equilibrium Tracking for \jpang{Prosumer} Energy Trading Games}
	
	\author{Yongkai Xie, Zhaojian Wang, John Z.F. Pang, Bo Yang, and Xinping Guan% <-this % stops a space
		\thanks{This work was supported by the National Natural Science Foundation of China (No. 62103265), and the ``Chenguang Program" supported by the Shanghai Education Development Foundation and Shanghai Municipal Education Commission of China (20CG11). (\textit{Corresponding author: Zhaojian Wang})	}	% <-this % stops a space
		\thanks{Y. Xie, Z. Wang, B. Yang, and X. Guan are with the Key Laboratory of System Control, and Information Processing, Ministry of Education of China, Department of Automation, Shanghai Jiao Tong University, Shanghai 200240, China, (email:\tt\small wangzhaojian@sjtu.edu.cn). }%%
		\thanks{J.Z.F. Pang is with the Institute of High Performance Computing (IHPC), A*STAR, Singapore 138632, Singapore, (email: \tt\small john\_pang@ihpc.a-star.edu.sg). }
	}
	
	\maketitle
	
	\begin{abstract}
		With the proliferation of distributed generations, traditional passive consumers in distribution networks are evolving into ``prosumers", which can both produce and consume energy. Energy trading with the main grid or between prosumers is inevitable if the energy surplus and shortage exist. To this end, this paper investigates the peer-to-peer (P2P) energy trading market, which is formulated as a generalized Nash game. We first prove the existence and uniqueness of the generalized Nash equilibrium (GNE). Then, an distributed online algorithm is proposed to track the GNE in the time-varying environment. Its regret is proved to be bounded by a sublinear function of learning time, which indicates that the online algorithm has an acceptable accuracy in practice. Finally, numerical results with six microgrids validate the performance of the algorithm.
	\end{abstract}
	
	
	\section{Introduction}
	
	% 	\subsection{Background}
	The explosive growth of distributed \jpang{generation} in distribution networks together with the advancement of communication and control technology at the consumer level have gradually transformed \jpang{the traditionally} passive consumers into ``prosumers", which can both produce and consume energy \cite{chen2019energy}. 
	Then, energy trading with the main grid or between prosumers is inevitable \jpang{since energy surplus and shortage are bound to exist} \cite{Wang2021Distributed}. In this situation, the peer-to-peer (P2P) market, \jpang{which operates} in a distributed manner, is more popular due to the ever-increasing number of prosumers, \jpang{in which} the various prosumers can be self-organized to operate economically and reliably under a given market mechanism \cite{Morstyn2018peer-to-peer}. In addition, the \jpang{increasing penetration and an aggravating volatility of} renewable generation calls for online market clearing methods. In this paper, we intend to investigate the distributed online energy trading market \jpang{for prosumers}. 
	
	% 	\subsection{Related Work}
	For such P2P energy trading markets, they are usually formulated as generalized Nash games, where each prosumer maximizes its profit with coupling constraints, \jpang{e.g.,} global power balance \cite{chen2019energy,Wang2021Distributed,Cui2020FairPeer-to-Peer,Nespoli2018rational,Chen2021SameFlexibility,Belgioioso2022Operationally}. Then, clearing the resulting P2P market corresponds to finding the generalized Nash equilibrium (GNE) of the energy trading game. For example, in \cite{chen2019energy}, the energy sharing game among prosumers is formulated with full information, and \cite{Wang2021Distributed} further designs a fully distributed algorithm based on Nesterov's methods to seek the GNE with only partial-decision information. In \cite{Cui2020FairPeer-to-Peer}, a P2P energy market is formulated as a generalized Nash game, where the \jpang{prosumers who share payments} are mutually coupled and influenced. \jpang{Following this,} \cite{Nespoli2018rational} and \cite{Chen2021SameFlexibility} further consider system-level grid constraints. \jpang{Lastly, in} \cite{Belgioioso2022Operationally}, a P2P energy market of prosumers is formulated as a generalized aggregative game with global coupling constraints.
	The aforementioned works have made great progress in the distributed GNE seeking for the P2P energy trading market. However, they usually focus on \jpang{only} one time section and provide offline solutions to solve the game. Due to the volatility of renewable generations and the complexity of load profiles, \jpang{both current and future operation status changes much more over time, requiring much faster algorithms}, i.e., online GNE tracking.
	
	% 	\subsection{Contributions}
	In this paper, we formulate a P2P energy trading market among prosumers in the distribution network and propose a distributed online algorithm to track the GNE of the market. The major contributions are as follows.
	\begin{itemize}
		\item A P2P energy trading market is modeled as a generalized Nash game with both individual and \jpang{coupled} time-varying constraints. Moreover, we prove the uniqueness of the GNE of this market at any time section. 
		\item A novel distributed online algorithm is proposed to track the GNE, where each prosumer can make decisions only using local variables and neighboring information. This reduces the communication burden and makes it easier to implement in practice. 
		\item \jpang{We prove a sublinear regret bound, i.e.,} that the regret of the online algorithm can be bounded by a sublinear function of learning time, indicating that the online algorithm suffers minimal ``loss in hindsight''. 
	\end{itemize}
	
	% 	\subsection{Structure}
	The rest of this paper is organized as follows. In Section \uppercase\expandafter{\romannumeral2}, the P2P energy trading game is formulated. \jpang{Section \uppercase\expandafter{\romannumeral3} introduces and analyzes the performance} of a distributed online algorithm to track the GNE of the game in a time-varying environment. 
	% 	The algorithm performance is also analyzed in Section \uppercase\expandafter{\romannumeral3}.
	Numerical results are presented in Section \uppercase\expandafter{\romannumeral4} to verify the effectiveness of our algorithm. Finally, Section \uppercase\expandafter{\romannumeral5} concludes the paper.
	
	\textit{Notations}: 
	In this paper, $\mathbb{R}_{+}^n$ is the $n$-dimensional (nonpositive) Euclidean space. For a column vector $x\in \mathbb{R}^n$ (matrix $A_{m\times n}\in \mathbb{R}^{m\times n}$), its transpose is denoted by $x^\mathsf{T}$($A^\mathsf{T}$). For a matrix $A$, $ [A]_{i,j}$ stands for the entry in the $i$-th row and $j$-th column of A.
	For vectors $x,y\in \mathbb{R}^n$, $x^\mathsf{T}y=\left\langle x,y \right\rangle$ is the inner product of $x,y$, while $\otimes$ represents the Kronecker product. 
	$\left\|x \right\|=\sqrt{x^\mathsf{T}x}$ is the Euclidean norm. 
	The identity matrix with dimension $n$ is denoted by $I_n$. Sometimes, we also omit $n$ to represent the identity matrix with the proper dimension. $\mathbf{0}_{n}, \mathbf{1}_{n}$ are all zero and all one vectors with dimension $n$, respectively. The Cartesian product of the sets $\Omega_i, i=1, \cdots, n$ is denoted by $\prod_{i=1}^n\Omega_i$. 
	Given a collection of $y_i$ for $i$ in a certain set $Y$, the vector composed of $y_i$ is defined as ${y}={col}(y_{i}):=(y_1, y_2, \cdots, y_n)^\mathsf{T}$. The projection of $x$ onto a set $\Omega$ is defined as $  	\mathcal{P}_{\Omega}(x):=\arg \min_{y\in \Omega}\left\|x-y \right\|$.

	\section{Problem Formulation}
	% 	In this section
	\subsection{Network model}
	
	We consider a distribution network with a group of prosumers, denoted by the set $\mathcal{N}=\{1,2,...,N\}$. For each prosumer, its load demand can be satisfied by its own generation and trading with the main grid or its neighboring prosumers. The trading edge is denoted by $\mathcal{E}\subseteq \mathcal{N}\times \mathcal{N}$. For a prosumer $ i $, the set of its neighbors is denoted by $ \mathcal{N}_i =\{\mathcal{N}_1^1,\dots,\mathcal{N}_1^{N_i}\}$ with $\left|\mathcal{N}_i\right|=N_i$. If $ j\in \mathcal{N}_i $, prosumers $ i $ and $ j $ can trade and communicate directly. Otherwise, direct trading and communication are not allowed. Then, the trading network is modeled as an undirected graph $\mathcal{G}=(\mathcal{N},\mathcal{E})$. The adjacency matrix of $\mathcal{G}$ is denoted by $W$ with elements $w_{i,j}$. If $j\in \mathcal{N}_i$, the weight $w_{i,j}$ satisfies $w_{i,j}=w_{j,i}>0$. Otherwise, $w_{i,j}=w_{j,i}=0$.  The Laplacian matrix of the communication graph is denoted by $L$ and we have $\textbf{1}^{\rm T}L=0$, where $ \textbf{1} $ is an all-one vector. 
	Moreover, the graph $\mathcal{G}$ is assumed to be connected.
	% \begin{assumption}\label{Assumption_connectivity}
		%     $\mathcal{G}$ is connected.
		% \end{assumption}
	For the weights, we have the following assumption, which implies that every row sum of $W$ is identical.
	\begin{assumption}\label{assumption_weight}
		The weight $w_{i,i}>0$ and $\sum_{j\in \mathcal{N}}w_{i,j}=w_0>0$ for all $i\in \mathcal{N}$.
	\end{assumption}
	
	\subsection{Prosumer model }	
	The scenario is that each prosumer is equipped with dispatchable generation, a non-dispatchable load, and an energy storage system (ESS). To maintain power balance, it can generate electricity, \jpang{charge or discharge from} the ESS, \jpang{and/or} trade with the main grid or neighboring prosumers. In this paper, we focus on the time horizon $\mathcal{T}=\{1,2,...,T\}$. Here, we will introduce them in detail.
	
	\emph{Dispatchable generation:} The power generated by dispatchable generation units of prosumer $i$ at time $t$, denoted by $p_i^g(t)$, is limited by
	\begin{align}\label{generation_limit}
		p_i^{g,min}\leq p_i^g(t)\leq p_i^{g,max},\quad \forall i\in \mathcal{N},t \in \mathcal{T}
	\end{align}
	where $p_i^{g,min}$ and $p_i^{g,max}$ are minimum and maximum local generation, respectively. Its generation cost is as follows.
	\begin{align}\label{generation_cost}
		f_i^g\left(p_i^g(t)\right)=a_i^g \left(p_i^g(t)\right)^2+b_i^gp_i^g(t)
	\end{align}
	where $a_i^g>0$ and $b_i^g$ are constants.
	
	\emph{Energy Storage Systems (ESS):} The ESS profile is constrained by the following dynamics.
	\begin{align}
		&0\leq p_i^c(t)\leq p_i^{c,max},\quad \forall i\in \mathcal{N},t \in \mathcal{T}\label{charge_limit}\\
		&0\leq p_i^d(t)\leq p_i^{d,max},\quad \forall i\in \mathcal{N},t \in \mathcal{T}\label{discharge_limit}\\
		&s_i(t+1)=s_i(t)+\frac{\Delta t}{e_i^{cap}}\left(\eta_i^c p_i^c(t)-\frac{1}{\eta_i^d}p_i^d(t)\right)\label{soc}\\
		&s_i^{min}\leq s_i(t)\leq s_i^{max},\quad \forall t\in \mathcal{T}\label{soc_limit}
	\end{align}
	where $p_i^c(t)$, $p_i^d(t)$, and $s_i(t)$ are the charging, discharging power, and state of charge (SoC) of the ESS $i$, respectively. $\Delta t$, $e_i^{cap}$, $\eta_i^c$ and $\eta_i^d$ are sampling time, ESS maximum storage capacity, and (dis)charging efficiencies, respectively. Moreover, $s_i^{min}$ and $s_i^{max}$, with $0<s_i^{min}<s_i^{max}<1$, denote the minimum and maximum SoC, while $p_i^{c,max}$ and $p_i^{d,max}$ denote the maximum (dis)charging power.

	Each prosumer might also minimize the usage of its ESS to reduce its degradation. Depending on the efficiency of a storage unit, there are losses based on usage that usually grow quadratically in power. \jpang{For simplicity, we disregard the effects due to SoC levels.} As defined in \cite{Hans2019Hierarchical}, the corresponding cost function is
	\begin{align}\label{es_cost}
		f_i^{es}\left(p_i^c(t),p_i^d(t)\right)=a_i^c \left(p_i^c(t)\right)^2+a_i^d \left(p_i^d(t)\right)^2
	\end{align}
	where $a_i^c$ and $a_i^d$ are both positive constants. 
	
	\emph{Trading with the main grid:}  Let $p_i^{mg}(t)$ be the power purchased from the main grid at time $t$ and $p^{mg}(t)=col\left\{p_i^{mg}(t)\right\}_{i\in\mathcal{N}}$. Similar to \cite{Atzeni2013Demand-Side}, \jpang{we set grid cost as}
	\begin{align}\label{grid_cost}
		C_t(p^{mg}(t))=c_t^{mg}\left(\sum\nolimits_{i\in \mathcal{N}}p_i^{mg}(t)\right)^2
	\end{align}
	where $c_t^{mg}$ is a \jpang{time-varying cost coefficient}, since the energy production varies along the time period according to the energy demand and the availability of distributed energy sources. Then the cost assigned to prosumer $i$ is
	\begin{align}\label{mg_cost}
		f_i^{mg}(p^{mg}(t))&=\frac{p_i^{mg}(t)}{\sum_{i\in \mathcal{N}}p_i^{mg}(t)}C_t(p^{mg}(t))\nonumber\\
		&=c_t^{mg}p_i^{mg}(t)\sum\nolimits_{i\in \mathcal{N}}p_i^{mg}(t)
	\end{align}
	Moreover, the total power exchanged with the main grid is limited, i.e.,
	\begin{align}\label{mg_coupling}
		p^{mg,min}\leq \sum\nolimits_{i\in \mathcal{N}}p_i^{mg}(t)\leq p^{mg,max},\quad \forall t \in \mathcal{T}
	\end{align}
	
	\emph{Trading with neighbors:} The trading cost with neighboring prosumers \jpang{of prosumer $i$} is 
	\begin{align}\label{trade_cost}
		f_i^{tr}\left(\left\{p_{i,j}^{tr}(t)\right\}\right)=\sum_{j\in \mathcal{N}_i}\left[a^{tr}\left(p_{i,j}^{tr}(t)\right)^2+d_{i,j}^{tr}p_{i,j}^{tr}(t)\right]
	\end{align}
	where $p_{i,j}^{tr}(t)$ is the power purchased from prosumer $j$ at time $t$, $d_{i,j}^{tr}=d_{j,i}^{tr}>0$ is the price and $a^{tr}$ is a small positive constant, which represents the tax cost incurred by using the energy sharing platform. 
	
	Disregarding loss on the power lines, the sum of the trading power of prosumer $i$ and $j$ at time $t$ should be 0.
	\begin{align}\label{trade_coupling}
		p_{i,j}^{tr}(t)+p_{j,i}^{tr}(t)=0,\quad \forall (i,j)\in \mathcal{E},t \in \mathcal{T}
	\end{align}
	Furthermore, trade between prosumers is limited by
	\begin{align}\label{trade_limit}
		p_{i,j}^{tr,min}\leq p_{i,j}^{tr}(t)\leq p_{i,j}^{tr,max},\quad \forall (i,j)\in \mathcal{E},t \in \mathcal{T}
	\end{align}
	where $p_{i,j}^{tr,min}\leq 0$ and $p_{i,j}^{tr,max}\geq 0$ are the minimum and maximum \jpang{tradeable power between prosumers $i$ and $j$}.
	
	Denote by $p_i^{l}(t)$ the undispatchable load demand, and the local power balance for each prosumer $i$ is
	\begin{align}\label{power_balance}
		&p_i^g(t)-p_i^c(t)+p_i^d(t)+p_i^{mg}(t)+\sum\nolimits_{j\in \mathcal{N}_i}p_{i,j}^{tr}(t)=p_i^{l}(t)
	\end{align}
	\subsection{Energy trading game}
	Before giving the game model, we first simplify notations. The decision variable of prosumer $i$ is denoted by
	\begin{align*}
		x_i(t):=col\left(p_i^g(t),p_i^c(t),p_i^d(t),p_i^{mg}(t),{p_{i,j}^{tr}(t)}\right)\in\mathbb{R}^{n_i}
	\end{align*}
	where $n_i=4+N_i$ and $\sum\nolimits_{i\in\mathcal{N}}n_i=n$. 
	
	Define a sparse matrix $E_i$, where the rows of $E_i$ correspond to every trading edge in $\mathcal{E}$ one by one. Let the $k$-th row of $E_i$ corresponds to $\left(I_k,J_k\right)$ in $\mathcal{E}$, then the elements $[E_i]_{k,l}$ of $E_i$ are assigned as follows.
	\begin{align*}
		[E_i]_{k,l}=\left\{
		\begin{matrix}
			1,&\text{If } \left\{I_k,J_k\right\}=\left\{i,\mathcal{N}_i^l\right\} \text{ and } I_k<J_k. \\
			-1,&\text{If } \left\{I_k,J_k\right\}=\left\{i,\mathcal{N}_i^l\right\} \text{ and } I_k>J_k.\\
			0,&\text{Otherwise.}
		\end{matrix}
		\right.
	\end{align*}
	
	Let
	\begin{align*}
		A_i&=\begin{bmatrix}
			\mathbf{0}_{m\times 3}&
			\begin{matrix}
				-1&\mathbf{0}^\mathsf{T}_{N_i}\\
				1&\mathbf{0}^\mathsf{T}_{N_i}\\
				\mathbf{0}_{\sum_{i\in \mathcal{N}}N_i}&E_i
			\end{matrix}
		\end{bmatrix}\\
		b_i&=\begin{bmatrix}
			-\frac{p^{mg,min}}{N}& 
			\frac{p^{mg,min}}{N}& 
			\mathbf{0}^\mathsf{T}_{\sum_{i\in \mathcal{N}}N_i}
		\end{bmatrix}^\mathsf{T}\\
		g_i(x_i(t))&=A_ix_i(t)-b_i,\quad \forall i\in \mathcal{N},t \in \mathcal{T}
	\end{align*}
	where $m=2+\sum_{i\in \mathcal{N}}N_i$. 
	Then coupling constraints \eqref{mg_coupling} and \eqref{trade_coupling} can be reformulated as
	\begin{align}\label{coupling_constraint}
		\sum\nolimits_{j\in \mathcal{N}}g_i(x_i(t))\leq 0
	\end{align}
	
	In \eqref{coupling_constraint}, the sparse matrix $E_i$ is designed to address the equality constraint \eqref{trade_coupling} by transforming it into two equivalent inequality constraints.
	% 	\begin{align*}
		% 		p_{i,j}^{tr}(t)+p_{j,i}^{tr}(t)\leq0,\quad \forall (i,j)\in \mathcal{E},t \in \mathcal{T}\nonumber\\
		% 		-p_{i,j}^{tr}(t)-p_{j,i}^{tr}(t)\leq0,\quad \forall (i,j)\in \mathcal{E},t \in \mathcal{T}
		% 	\end{align*}
	% 	while elements $-1$, $1$ in $A_i$ and the previous two elements in $b_i$ correspond to constraint \eqref{mg_coupling}.
	
	Similarly, let 
	\begin{align*}
		G_i=\begin{bmatrix}
			1&-1&1&1&\mathbf{1}^\mathsf{T}_{N_i}
		\end{bmatrix}
	\end{align*}
	Then equality constraint $\eqref{power_balance}$ can be reformulated as 
	\begin{align}\label{equality_constraint}
		G_ix_i(t)-p_i^{l}(t)=0,\quad \forall i\in \mathcal{N},t \in \mathcal{T}
	\end{align}

	Moreover, the domain set of $x_i(t)$ is denoted by
	\begin{align*}
		&\chi_i=\left\{x_i(t)|x_i(t)\ \text{satisfies}\  \eqref{generation_limit},\eqref{charge_limit},\eqref{discharge_limit},\eqref{soc},\eqref{soc_limit},\eqref{trade_limit}\right\}\\
		&\Omega_i=\chi_i\cap \left\{x_i(t)|x_i(t)\text{ satisfies }\eqref{equality_constraint}\right\}\\
		&X=\prod\nolimits_{i\in\mathcal{N}}\Omega_i\cap \left\{x(t)|x(t)\text{ satisfies }\eqref{coupling_constraint}\right\}
	\end{align*}
	where $\chi_i$ is the set of all of the local inequality constraints, $\Omega_i$ considers the \jpang{reformulated power balance constraint of \eqref{power_balance}}, while $X$ includes coupling constraints. 
	
	In the energy trading game, each prosumer intends to minimize its cost while maintaining the global power balance. The optimization problem of each prosumer is 
	\begin{subequations}\label{Energy_trading_game}
		\begin{align}\label{Energy_trading_game_obj}
			\min_{x_i(t)} \ &J_{i,t}\left(x_i(t),x_{-i}(t)\right)=f_i^g\left(p_i^g(t)\right)+f_i^{es}\left(p_i^c(t),p_i^d(t)\right)\nonumber\\
			&\qquad+f_i^{mg}\left(p^{mg}(t)\right)+f_i^{tr}\left(\left\{p_{i,j}^{tr}(t)\right\}\right)\\
			\text{ s.t.}\ & x_i(t)\in X_i(x_{-i}(t))
		\end{align}
	\end{subequations}
	where $x_{-i}(t):=\operatorname{col}\left(x_1(t),...,x_{i-1}(t),x_{i+1}(t),...,x_N(t)\right)$ and $X_i(x_{-i}(t)):=\left\{x_i(t)\in \Omega_i| (x_i(t),x_{-i}(t))\in X\right\}$. 
	
	In summary, the energy trading game is represented as  
	\begin{itemize}
		\item Player: all prosumers, denoted by $\mathcal{N}=\{1,2,...,N\}$.
		\item Strategy: decision variable $x_i$.
		\item Payoff: the disutility function $J_{i,t}\left(x_i(t),x_{-i}(t)\right)$.
	\end{itemize}
	Due to the global coupling constraints \eqref{coupling_constraint}, it is a generalized Nash game. The corresponding GNE is defined as
	\begin{align*}
		x_i^*(t)&\in \arg \min\; J_{i,t}\left(x_i(t),x_{-i}^*(t)\right),  \text{s.t. } x_i(t)\in X_i\left(x_{-i}^*(t)\right)
	\end{align*}
	
	Regarding the game \eqref{Energy_trading_game}, we have following assumptions.
	\begin{assumption}\label{assumption_domain_set}
		$\chi_i$ is a  non-empty, compact and convex set.
	\end{assumption}
	\begin{assumption}\label{assumption_feasibility}
		Given any $x_{-i}(t)$, problem \eqref{Energy_trading_game} is feasible.
	\end{assumption}
	Since the constraints of problem \eqref{Energy_trading_game} are all affine, the commonly assumed Slater's condition is simplified to only require feasibility, and therefore Assumption \ref{assumption_feasibility} suffice. %Thus, Assumption \ref{assumption_feasibility} is a common assumption, otherwise, the problem is infeasible. 

	\subsection{Uniqueness of the GNE}
	The pseudo-gradient of $\{J_{i,t}\}_{i\in\mathcal{N}}$ is defined as
	\begin{align}\label{pseudo_gradient_definition}
		&F_t(x)=col\left(\left\{\nabla_{x_i}J_{i,t}(x_i,x_{-i})\right\}\right)\nonumber\\
		&=col\left(\left\{\begin{matrix}
			2a_i^gp_i^g+b_i^g\\
			2a_i^cp_i^c\\
			2a_i^dp_i^d\\
			c_t^{mg}\left(2p_i^{mg}+\sum_{j\in \mathcal{N},j\neq i}p_j^{mg}\right)\\
			col\left(\left\{2a^{tr}p_{i,j}^{tr}+d_{i,j}^{tr}\right\}_{j\in \mathcal{N}_i}\right)
		\end{matrix}\right\}\right)
	\end{align}
	Define $\underline {a}^g=\underset{i\in\mathcal{N}}{\min}\,a_i^g$, $\overline {a}^g=\underset{i\in\mathcal{N}}{\max}\,a_i^g$, \jpang{with $\underline {a}^c$, $\overline {a}^c$, $\underline {a}^d$, $\overline {a}^d$ defined similarly}. In addition, $\underline {c}^{mg}=\underset{t\in \mathcal{T}}{\min}\, c_t^{mg}$, $\overline {c}^{mg}=\underset{t\in \mathcal{T}}{\max}\, c_t^{mg}$. Then, we can prove that the pseudo-gradient is strongly monotone and Lipschitz continuous. 
	\begin{lemma}\label{lemma_monotonicity_continuity}
		For $\forall t\in \mathcal{T}$, the pseudo-gradient $F_t(x)$ has following properties\jpang{:}
		\begin{enumerate}[]
			\item $F_t(x)$ is strongly monotone with parameter $0<\eta\leq \min\left\{2\underline {a}^g,2\underline {a}^c,2\underline {a}^d, \underline {c}^{mg},2a^{tr}\right\}$, i.e., $\langle F_t(x)-F_t(y),x-y\rangle\geq \eta||x-y||_2^2$;
			\item  $F_t(x)$ is $\theta-$Lipschitz continuous, i.e., $\left\|F_t(x)-F_t(y)\right\|_2\leq \theta ||x-y||_2$ with $\theta\geq \max\left\{2\overline {a}^g,2\overline {a}^c,2\overline {a}^d,N\overline {c}^{mg},2a^{tr}\right\}$.
		\end{enumerate}
	\end{lemma}
	
	\begin{proof} 
		For 1), taking any two variables $x^1, x^2$, then
		\begin{align}
			&\langle F_t(x^1)-F_t(x^2),x^1-x^2\rangle\nonumber\\
			&=\sum\nolimits_{i \in \mathcal{N}}2a_i^g\left\|p_i^{g,1}-p_i^{g,2}\right\|_2^2+\sum\nolimits_{i \in \mathcal{N}}2a_i^c\left\|p_i^{c,1}-p_i^{c,2}\right\|_2^2\nonumber\\
			&\quad+\sum\nolimits_{i \in \mathcal{N}}2a_i^d\left\|p_i^{d,1}-p_i^{d,2}\right\|_2^2+c_t^{mg}h^THh\nonumber\\
			&\quad+\sum\nolimits_{i \in \mathcal{N}}\sum\nolimits_{j \in \mathcal{N}_i}2a^{tr}\left\|p_{i,j}^{tr,1}-p_{i,j}^{tr,2}\right\|_2^2
		\end{align}
		where $h=col\left(p_i^{mg,1}-p_i^{mg,2}\right)$ and $H=I_N+\textbf{1}_{N\times N}$. Since the minimal eigenvalue of $H$ is $1$, we have $h^THh\geq \sum_{i \in \mathcal{N}}\left\|p_i^{mg,1}-p_i^{mg,2}\right\|_2^2$. Then $\langle F_t(x^1)-F_t(x^2),x^1-x^2\rangle\geq\eta\left\|x^1-x^2\right\|_2^2$ with $0<\eta\leq \min\left\{2\underline {a}^g,2\underline {a}^c,2\underline {a}^d, \underline {c}^{mg},2a^{tr}\right\}$. Therefore, $F_t(x)$ is strongly monotone with parameter $\eta$.
		
		Similarly, the second assertion could be obtained, which is omitted due to the space limit. 
	\end{proof} 
	A GNE with the same Lagrangian multipliers for all the agents is called variational GNE (v-GNE) \cite{Facchinei2010Generalized}, which has the economic interpretation of no price discrimination \cite{Kulkarni2012variational}. For $\forall t\in\mathcal{T}$,  every solution $x^*(t)\in X$ to the following variational inequality is a v-GNE of game $\eqref{Energy_trading_game}$.
	\begin{align}\label{variational _inequality}
		\langle F_t(x^*(t)),x-x^*(t)\rangle\leq 0,\qquad \forall x\in X
	\end{align}
	
	% 	For the existence and uniqueness of the v-GNE, we have the following results.
	\begin{theorem}\label{lemma_equilibrium_existence}
		For the generalized Nash game $\eqref{Energy_trading_game}$, there exists a unique v-GNE.
	\end{theorem}
	\begin{proof}
		Following \cite{Facc2010Nashequilibria}, since $J_{i,t}(x_i(t),x_{-i}(t))$ is differentiable and convex with respect to $x_i(t)$ for any $x_{-i}(t)$, if {Assumption} {\ref{assumption_domain_set}} and {\ref{assumption_feasibility}} hold, the v-GNE of $\eqref{Energy_trading_game}$ exists. Moreover, by the strong monotonicity of $F_t(x)$, the uniqueness of v-GNE is guaranteed.
	\end{proof}
	
	\section{Online Algorithm}
	
	In this section, we first propose an online distributed algorithm based on a consensus algorithm and primal-dual strategy to solve the problem \eqref{Energy_trading_game}. Then, we prove that the regret of the proposed algorithm is bounded by a sublinear function of the learning time.
	\subsection{Algorithm design}
	Recalling the objective function \eqref{Energy_trading_game_obj}, $f_i^{mg}\left(p^{mg}(t)\right)$ is associated with decisions of all of the other prosumers. To solve game \eqref{Energy_trading_game}, full information is needed, i.e., one prosumer needs to communicate with all of the other prosumers. However, this is difficult to realize in practice due to communication limits. This section designs a distributed algorithm with only partial-decision information, where the prosumers only need to exchange information with its neighbors. To this end, we endow each prosumer with an auxiliary variable $\bar{x}^i(t)$ that provides an estimate of the decisions of other prosumers at time $t$. Moreover, \jpang{$\bar{x}_j^i(t)$ represents prosumer $i$'s estimate of $j$'s decision} and $\bar{x}^i_{-i}(t)=col\left(\left\{\bar{x}_j^i(t)\right\}_{j\in\mathcal{N}_i}\right)$.  Clearly, we have $\bar{x}_i^i(t)=x_i(t)$. 
	
	\jpang{Firstly, note that the} Lagrangian of problem \eqref{Energy_trading_game} is
	\begin{align}\label{Lagrangian_function}
		&L_{i,t}\left(x_i(t),\lambda_i(t),\mu_i(t);x_{-i}(t)\right)=J_{i,t}\left(x_i(t),x_{-i}(t)\right)\nonumber\\
		&\ +\lambda_i^T(t)\sum\nolimits_{i \in \mathcal{N}}g_i(x_i(t))+\mu_i(t)(G_ix_i(t)-p_i^l(t))
	\end{align}
	where $\lambda_i(t)$ and $\mu_i(t)$ are Lagrange multipliers.
	
	The iterative process of $x_i(t)$, $\bar{x}_{-i}^i(t)$ and dual variables $\left(\lambda_i(t),\mu_i(t)\right)$ is shown in Algorithm \ref{alg1}, where $0<\rho(t)<1$ is a stepsize or the so-called learning rate, which decreases over $t$, and  $c=\frac{1}{w_0}$.  
	\begin{algorithm}[t]
		\caption{Online Algorithm for P2P Energy Trading} 
		\label{alg1} 
		\begin{algorithmic} 
			\State \textbf{Initialization:} $x_i(0)\in \chi_i$, $\bar{x}_{-i}^i(0)\in\mathbb{R}^{n-n_i}$, $\lambda_i(0)=0$,
			$\mu_i(0)=0$
			\For{$t=1$ to $T$}
			% 			\State 
			\begin{subequations}
				\begin{align}
					&x_i(t+1)=(1-\rho(t))x_i(t)\nonumber\\
					&\qquad+\rho(t)\mathcal{P}_{\chi_i}\left\{x_i(t)
					-\rho(t)[\nabla_{x_i}J_{i,t}(x_i(t),\bar{x}_{-i}^i(t))\right.\nonumber\\
					&\qquad+ \rho(t)(A_i^T\lambda_i(t)
					+G_i^T\mu_i(t))\nonumber\\
					&\qquad\left.+c\sum\nolimits_{j\in \mathcal{N}_i}(x_i(t)-\bar{x}_{i}^j(t))]\right\}\label{algorithm_xi}\\
					&\bar{x}_{-i}^i(t+1)=\bar{x}_{-i}^i(t)\nonumber\\
					&\qquad-c \rho(t) \sum\nolimits_{j \in \mathcal{N}_i}w_{i,j}\left(\bar{x}_{-i}^i(t)-\bar{x}_{-i}^j(t)\right)\label{algorithm_x-i}\\
					&\lambda_i(t+1)=\mathcal{P}_{\mathbb{R}_+^m}\bigg\{\left(1-\rho(t)\right) \sum\nolimits_{j \in \mathcal{N}_i\cup \{i\}} \frac{w_{i j}}{w_0} \lambda_j(t)\nonumber\\
					&\qquad+\rho(t)\left[A_i\left(2 x_i(t+1)-x_i(t)\right)-b_i\right] \bigg\}\label{algorithm_lambda}\\
					&\mu_i(t+1)=\left(1-\rho(t)\right) \mu_i(t)\nonumber\\
					&\qquad+\rho(t)\left[G_i\left(2 x_i(t+1)-x_i(t)\right)-p_i^l(t)\right]\label{algorithm_mu}
				\end{align}
			\end{subequations}
			\EndFor
		\end{algorithmic} 
	\end{algorithm}
	
	The update for $x_i(t)$ in Algorithm \ref{alg1} employs the projected primal-dual gradient decomposition method combined with the consensus approach \cite{wang2019distributed2,wang2019distributed,wang2019distributed_variation}. The update of $\bar{x}_{-i}^i(t)$ can be regarded as a discrete-time integration for the consensus error of the local estimation \cite{pavel2019distributed}. At each sampling time $t$, prosumer $i$ gets $\nabla_{x_i}J_{i,t}\left(x_i(t),x_{-i}(t)\right)$ by using $x_i(t), \bar{x}_{-i}^i(t)$ and updates $x_i(t+1)$ with \eqref{algorithm_xi}. Meanwhile, the estimation $\bar{x}_{-i}^i(t+1)$ is updated by communication with neighboring prosumers by \eqref{algorithm_x-i}. Then, dual variables $\left(\lambda_i(t),\mu_i(t)\right)$ are updated using the latest updated local information $x_i(t+1)$ with \eqref{algorithm_lambda} and \eqref{algorithm_mu}, respectively. %In addition, to update $\mu_i(t)$, the undispatchable load demand at the next sampling time is needed, which can be obtained through short-term load forecasting.\wang{The load forecast is unreasonable.}
	
	\begin{remark}
		Algorithm \ref{alg1} is fully distributed with only partial-decision information. Each prosumer makes a decision only based on local information and communication with its neighbors, which is easy to implement in practice.
		Compared with the existing work \cite{Kaihong2021Environments} only considering the time-varying objective function, we further include the time-varying constraints. Moreover, compared with the algorithm in \cite{2021arXiv210506200M}, the update is much simpler without the need of solving an optimization problem at each iteration, which reduces the computation cost. 
	\end{remark}

	\subsection{Regret analysis}
	In this subsection, we will prove that the regret of Algorithm \ref{alg1} is bounded by a sublinear function of the learning time. First, we give some notations.
	Under {Assumption} {\ref{assumption_domain_set}},  $\left\|x_i\right\|$, $\left\|g_i(x_i)\right\|$, $\left\|\nabla_{x_i}J_{i,t}(x_i,x_{-i})\right\|$, $\left\|\nabla_{x_i}g_i(x_i)\right\|$ and  $\left\|G_ix_i-p_i^l(t)\right\|$ are bounded. Then, for $\forall i\in \mathcal{N}$ and $\forall t\in \mathcal{T}$, we define
	\begin{align*}
		\kappa_1&=\sup\nolimits_{x_i\in\chi_i}\left\|x_i\right\|,\ 			\kappa_2=\sup\nolimits_{x_i\in\chi_i}\left\|g_i(x_i)\right\|\\
		\kappa_3&=\sup\limits_{x_i\in\chi_i}\left\|\nabla_{x_i}J_{i,t}(x_i,x_{-i})\right\|,\
		\kappa_4=\sup\limits_{x_i\in\chi_i}\left\|\nabla_{x_i}g_i(x_i)\right\|\\
		\kappa_5&=\sup\limits_{x_i\in\chi_i}\left\|G_ix_i-p_i^l(t)\right\|,\
		\kappa_6=\max\limits_{i\in\mathcal{N}}\sup_{\left\|x_i\right\|=1}\left\|G_ix_i\right\|
	\end{align*}
	
	We use the dynamic regret to evaluate the performance of the online algorithm, which is defined as follows. 
	\begin{align}\label{regret_definition}
		R_i(T)=\sum\nolimits_{t=1}^T\left[J_{i,t}\left(x_i(t), x_{-i}^*(t)\right)-J_{i,t}\left(x^*(t)\right)\right] 
	\end{align}
	where $x^*(t)=col(x_i^*(t))$ is the v-GNE of \eqref{Energy_trading_game} at time $t$. 
	
	An online algorithm is \jpang{generally considered to perform well} if the regret increases sublinearly \cite{2021arXiv210506200M}, \cite{Shahrampour2018MirrorDescent}, i.e.,
	\begin{align}\label{regret_limit}
		\lim _{T \rightarrow \infty} \frac{R_i(T)}{T}=0,\quad \forall i \in \mathcal{N}
	\end{align}
	
	However, it would be impossible to keep the regret \eqref{regret_definition} increasing sublinearly if the v-GNE sequence $\left\{x^*(0),x^*(1),\dots,x^*(N)\right\}$ of \eqref{Energy_trading_game} fluctuates drastically. Therefore, motivated by \cite{Hall2015Online} and \cite{Kaihong2020Pseudoconvex}, we adopt the following accumulation to describe the difficulty of tracking the v-GNE sequence:
	\begin{align}\label{varphi_T_definition}
		\varPhi_T=\sum\nolimits_{t=0}^T\left\|x^*(t+1)-x^*(t)\right\|
	\end{align}
	Sublinear regret is only possible when $\varPhi_T$ is small.
	
	The next result shows that by implementing decreasing $\rho(t)$, the regret of Algorithm \ref{alg1} is bounded by $\varPhi_T$ and $\sum_{t=1}^T\sqrt{\rho(t)}$.
	\begin{theorem}\label{Theorem_1}
		Suppose that {Assumption} {\ref{assumption_weight}} and {\ref{assumption_domain_set}} hold, for $\forall i\in \mathcal{N}$, the regret of Algorithm \ref{alg1} is bounded by
		\begin{align}\label{theorem_1}
			R_i(T) \leq \mathcal{O}\left(\sqrt{T\left(\frac{\varPhi_T+1}{\rho^2(T)}+\sum\nolimits_{t=1}^T\sqrt{\rho(t)}\right)}\right)
		\end{align}
	\end{theorem}
	The proof of {Theorem \ref{Theorem_1}} is given in the appendix.
	\begin{remark}
		From {Theorem \ref{Theorem_1}}, we can get the sublinearly bounded dynamic regrets of Algorithm \ref{alg1}. To this end, take $\rho(t)=\frac{K}{(at+b)^\alpha}$ with $a,b>0$, $0<K\leq 1$ and $0<\alpha<\frac{1}{2}$. If $\varPhi_T$ is sublinear with $T^{1-\frac{\alpha}{2}}$, then $R_i(T)$ is sublinear with $T$. Note that the sublinearity of $\varPhi_T$ is a general assumption, which is widely adopted in online optimization and online game \cite{Kaihong2021Environments,2021arXiv210506200M,Kaihong2020Pseudoconvex}. In the P2P energy market, it mainly limits the fluctuation of the undispatchable load demand, which results in a GNE sequence with limited deviation. 
		% Thus, the regret of each prosumer is sublinear with $T$.
	\end{remark}

	% 	\begin{proof}
		% 		If $\rho(t)=\frac{K}{(at+b)^\alpha}$, then we have
		% % 		\begin{align}
			% % 			\sum_{t=1}^T\sqrt{\rho(t)}&<\int_{0}^{T}\sqrt{\frac{K}{(at+b)^\alpha}}dt\nonumber\\
			% % 			&\left.=\frac{\sqrt{K}}{a\left(1-\frac{\alpha}{2}\right)}\left(at+b\right)^{1-\frac{\alpha}{2}}\right|_0^T\leq \mathcal{O}\left(T^{1-\frac{\alpha}{2}}\right)\nonumber
			% % 		\end{align}
		% 		\begin{align*}
			% 			\sum_{t=1}^T\sqrt{\rho(t)}&<\int_{0}^{T}\sqrt{\frac{K}{(at+b)^\alpha}}dt\leq \mathcal{O}\left(T^{1-\frac{\alpha}{2}}\right)
			% 		\end{align*}
		
		% 		Thus, $R_i(T)\leq \mathcal{O}\left(T^{\frac{1}{2}+\alpha}\sqrt{\varPhi_T+1}+T^{1-\frac{\alpha}{4}}\right)$. Since $\varPhi_T$ is sublinear with $T^{1-\frac{\alpha}{2}}$, i.e., $\lim\limits_{T\rightarrow\infty}\frac{\varPhi_T}{T^{1-\frac{\alpha}{2}}}=0$, then $R_i(T)\leq \mathcal{O}(T)$, indicating that the regret of each prosumer increases sublinearly.
		% 	\end{proof}
	%\begin{remark}
	%\wang{Here, we need several sentences to explain the result of Theorem 2, such as the sublinearly bounded regret?}
	%\end{remark}
	
	\section{Case Studies}
	
	In this section, numerical simulation on a case with six prosumers (microgrids) is introduced to verify the effectiveness of the proposed algorithm. 
	%First,  is carried out to illustrate the major properties of our algorithm. 
	\begin{figure}[t]
		\centering
		\includegraphics[width=0.18\textwidth]{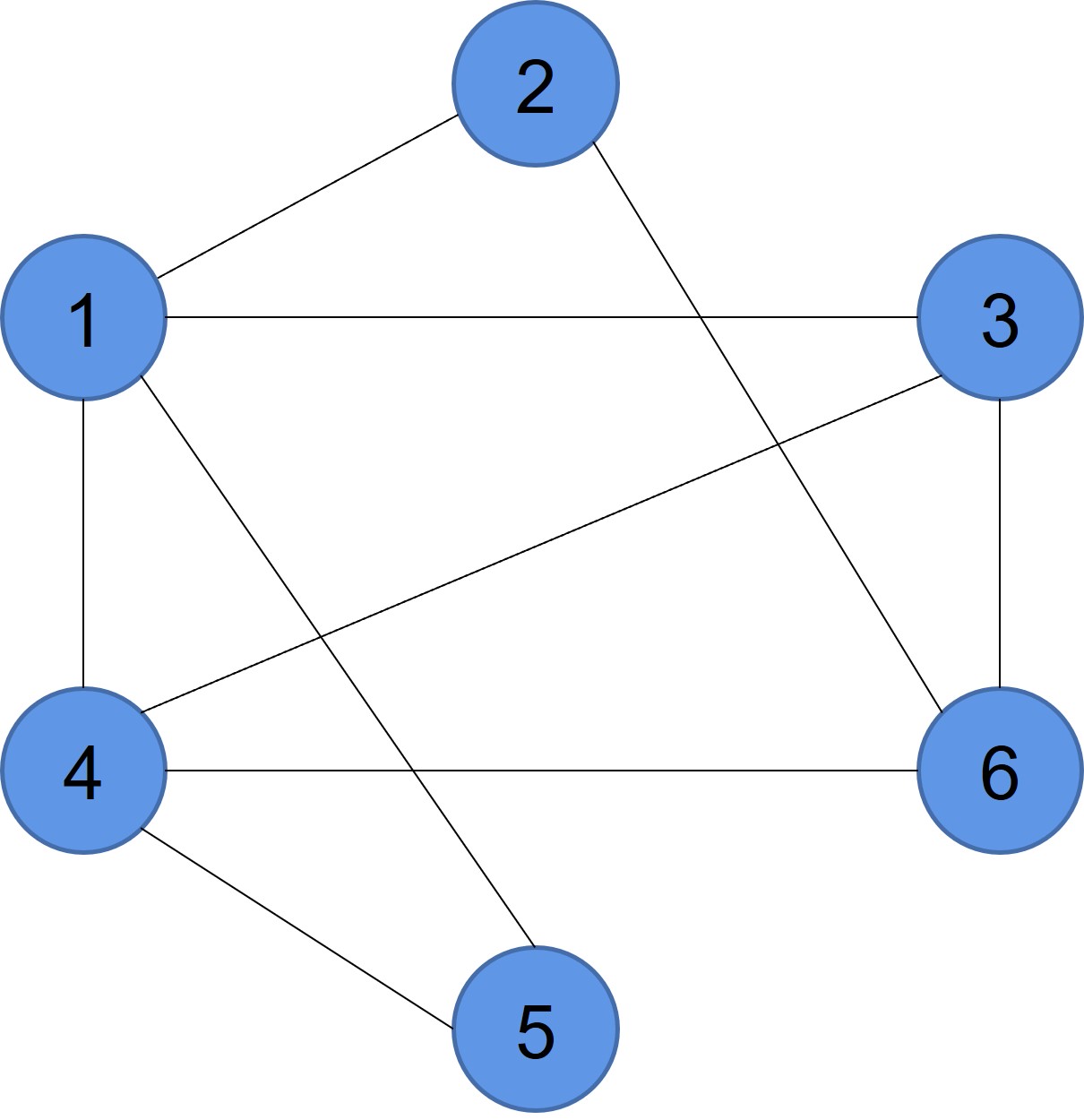}
		\caption{Communication and trading graph $\mathcal{G}$ for the case study.}
		\label{fig_topology}
	\end{figure}
	% 	\subsection{System Configuration} 
	Each prosumer communicates and trades with its neighboring prosumers via a connected graph $\mathcal{G}$, which is shown in Fig. \ref{fig_topology}. The adjacency matrix $W$ is set to be $w_{i,j}=w_{j,i}=\frac{1}{{N}_i+1}$ and $\sum_{j\in \mathcal{N}}w_{i,j}=1$ for $\forall i\in \mathcal{N}$. 
	The time interval is set to be $\Delta t=1\ \text{min}$ and $T=1440$, which indicates that the optimization period is 24 hours. The learning rate is set to be $0.8\sqrt[3]{\frac{5}{0.1t+5}}$.
	
\begin{figure}[t]
		\centering
		\includegraphics[width=0.3\textwidth]{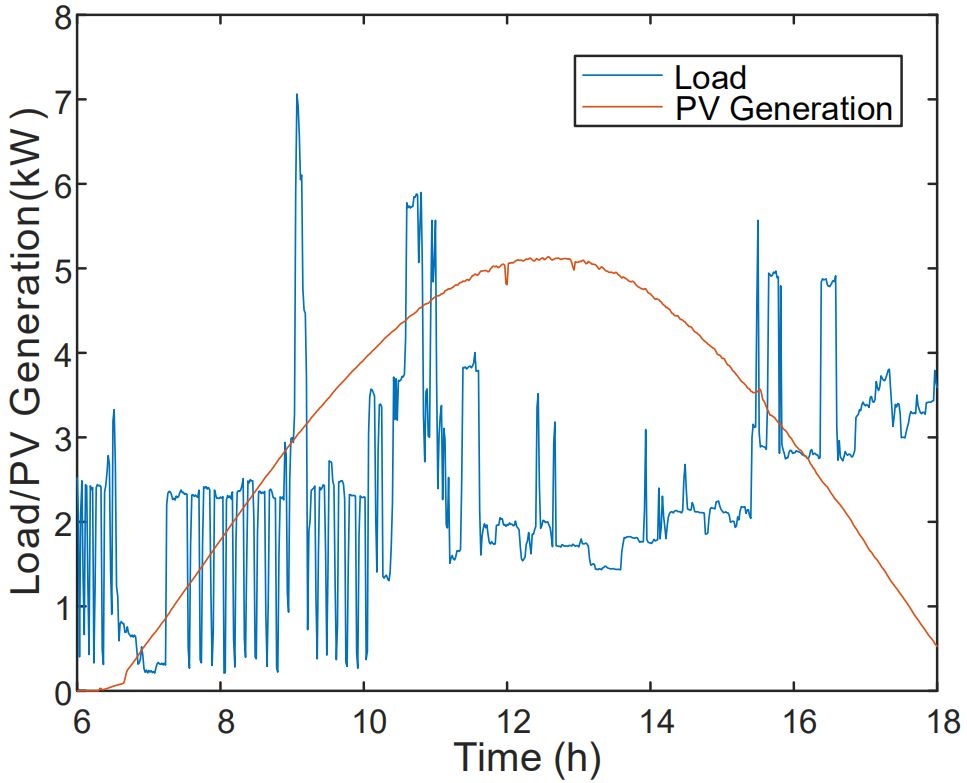}
		\caption{Load and PV generation of prosumer 1.}
		\label{fig_pv_load}
	\end{figure}

	We assume that prosumer 1, 3, 4, and 6 are equipped with PVs and $p_i^l(t)$ of these prosumers is equal to the difference between their loads and PV generations, while $p_i^l(t)$ of prosumer 2 and 5 only have a load demand. The minute-sampled profile of PV generation is obtained from \cite{website_PV}, which is collected in Utah, from the U.S.. We use the data on 16th September 2013. Since each prosumer is located in the same area, we assume that the photovoltaic generation curves of different prosumers are the same, but with different amplitudes. The profiles of the daily power consumption of each prosumer are from \cite{website_load}. Fig. $\ref{fig_pv_load}$ depicts the load and PV generation of prosumer 1 from 6:00 am to 6:00 pm.
	
	% 	\subsection{Simulation Results}
	\begin{figure}[t]
		\centering
		\includegraphics[width=0.45\textwidth]{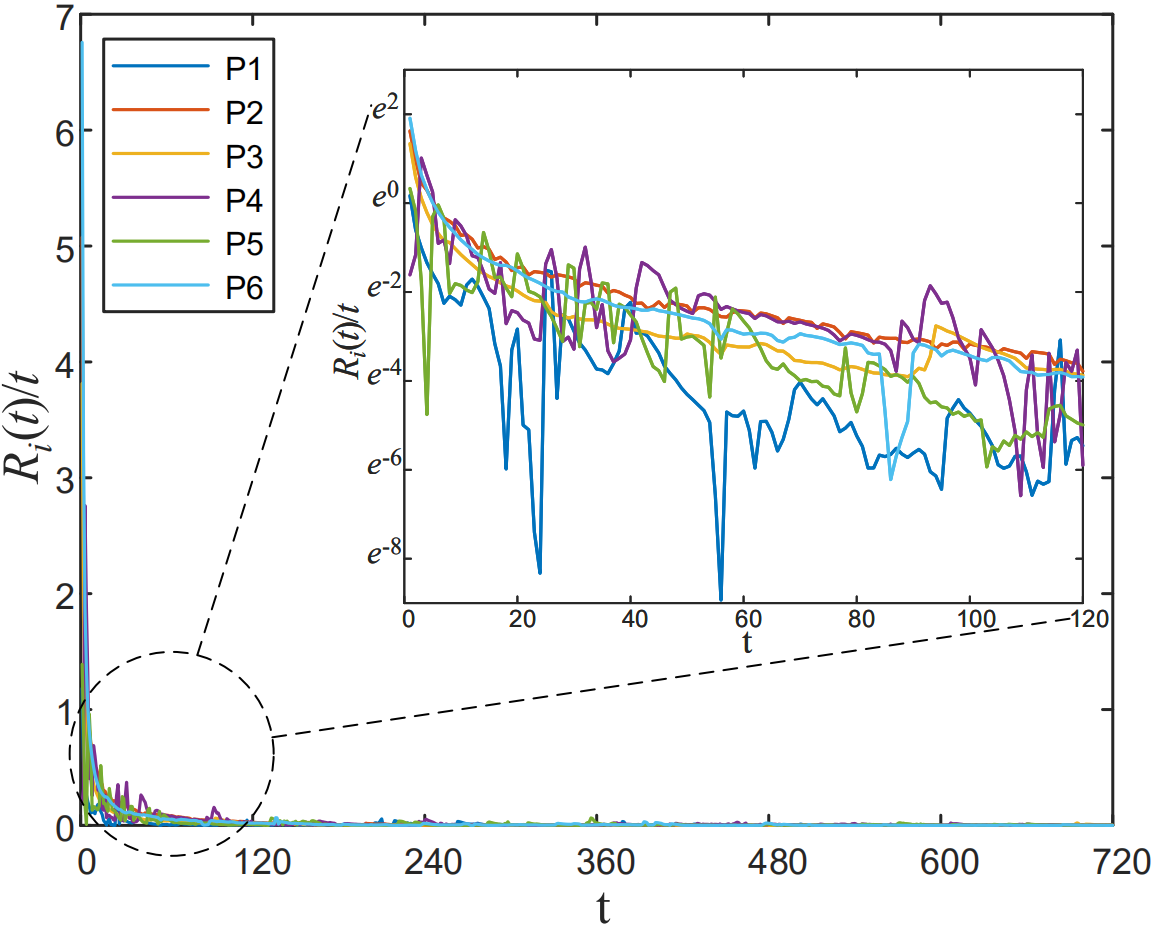}
		\caption{The trajectories of the average regrets of each prosumer. Note log-scale on y-axis for the inlay figure.}
		\label{fig_reg}
	\end{figure}	
	
	Fig. \ref{fig_reg} shows the trajectories of the average regrets of each prosumer from 6:00 am to 6:00 pm. We set $t=0$ at 6:00 am here and therefore the length of the learning period is 720. As shown in Fig. $\ref{fig_reg}$, $R_i(t)/t$, $i\in\mathcal{N}$ approximately decays to zero after $t=120$ (i.e., after two hours). The downtrend of the logarithm of $R_i(t)/t$, $i\in\mathcal{N}$ also validates this property, which is consistent with the result in {Theorem \ref{Theorem_1}}.

	% 	Fig. \ref{fig_mg_sum} shows the aggregated trading power with the main grid of each prosumer over time. As shown in Fig. $\ref{fig_mg_sum}$, $\sum_{i=1}^{6}p_i^{mg}(t)$ oscillates slightly due to the variation of the loads and PV generations, which still satisfies the upper and lower limits. 
	% 	\begin{figure}[t]
		% 		\centerline{\includegraphics[scale=0.35]{mg_sum_1.jpg}}
		% 		\caption{Aggregated trading power with the main grid of six prosumers.}
		% 		\label{fig_mg_sum}
		% 	\end{figure}

	% 	\subsection{Simulation on -node System}
	
	\section{Conclusion}
	In this paper, we propose a distributed online algorithm to track the GNE of the P2P energy trading game in a time-varying environment. We prove that by an appropriate choice of the decreasing learning rate, the regret of the proposed algorithm is bounded by a sublinear function of learning time. Simulation results with six prosumers verify the effectiveness of the proposed algorithm. In future work, we may focus on the effect of communication delay on the performance of online algorithms.

	\appendix[proof of theorem \ref{Theorem_1}]
	We start with a lemma.
	\begin{lemma}\label{lemma_varepsilon}
		If {Assumption \ref{assumption_weight}} holds,  we have $0<\varepsilon<1$, where $\varepsilon=\max \limits_{\substack{t\in\mathcal{T}, s_i>0}}\left\{\left|1-c \rho(t) s_i\right|\right\}$ with $0=s_1\leq s_2\leq \dots\leq s_N$ as $N$ eigenvalues of $L$.
	\end{lemma}
	
	\begin{proof}
		Let $d^*=\max\limits_{i\in\mathcal{N}}\left\{\sum_{j\in \mathcal{N}}w_{i,j}\right\}$. By {Assumption \ref{assumption_weight}}, $d^*=w_0$. From \cite[pp. 31]{Mesbahi2010Graph}, since $s_N$ is the maximal eigenvalue of $L$, we have $d^*\leq s_N\leq 2d^*$. Therefore, $w_0\leq s_N\leq 2w_0$. Recall $0<\rho(t)<1$ and $c=\frac{1}{w_0}$, we have $-1<1-c\rho(t)s_N<1$. Thus, for $\forall$ $t\in\mathcal{T}$ and $s_i>0$, we have $-1<1-c\rho(t)s_N\leq 1-c\rho(t)s_i<1$. Thus, $0<\varepsilon<1$.
	\end{proof}
	
	The estimation error of prosumer $i$ is defined as $e_i(t)=\operatorname{col}\left(e_i^1(t), \ldots, e_i^N(t)\right)$ with $e_i^j(t)=\bar{x}_i^j(t)-x_i(t)$. Based on {Lemma \ref{lemma_varepsilon}}, we now present the bound of $e_i(t)$. 
	\begin{lemma}\label{lemma_error_bound}
		Under {Assumption} {\ref{assumption_weight}} and {\ref{assumption_domain_set}}, for $\forall i\in \mathcal{N}$ and $2\leq t\leq T$, we have
		%		\begin{subequations}
			\begin{align*}
				||e_i(t)|| \leq \varepsilon^{t-1}||e_i(1)||+2 \sqrt{N} \kappa_1 \sum\nolimits_{k=0}^{t-2} \varepsilon^k \rho(t-k-1)
				%				\left\|e_i(t)\right\|^2 \leq \varepsilon^{t-1}\left\|e_i(1)\right\|^2+\frac{8 N \kappa_1^2}{1-\varepsilon} \sum_{k=0}^{t-2} \varepsilon^k \rho(t-k-1)
			\end{align*}
			%		\end{subequations}
	\end{lemma}
	\begin{proof}By \eqref{algorithm_x-i}, we have
		\begin{align}
			&e_i^j(t+1) =\bar{x}_i^j(t+1)-x_i(t+1) \nonumber\\
			&=\bar{x}_i^j(t)-c \rho(t) \sum\nolimits_{k \in \mathcal{N}_j}w_{i,j}\left(\bar{x}_i^j(t)-\bar{x}_i^k(t)\right)-x_i(t+1) \nonumber\\
			&=\bar{x}_i^j(t)-x_i(t)-\left(x_i(t+1)-x_i(t)\right) \nonumber\\
			& -c \rho(t) \sum\nolimits_{k \in \mathcal{N}_j}w_{i,j}\left[\left(\bar{x}_i^j(t)-x_i(t)\right)-\left(\bar{x}_i^k(t)-x_i(t)\right)\right]\nonumber\\
			&=e_i^j(t)-c \rho(t) \sum\nolimits_{k \in \mathcal{N}_j}w_{i,j}\left(e_i^j(t)-e_i^k(t)\right)\nonumber\\
			&-\left(x_i(t+1)-x_i(t)\right)
		\end{align}
		Let $\Delta x_i(t)=x_i(t+1)-x_i(t)$, then
		\begin{align*}
			e_i(t+1)=&e_i(t)-c\rho(t)\left(L \otimes I\right) e_i(t)-\left(\mathbf{1}_N\otimes I\right)\Delta x_i(t)
		\end{align*}
		By \eqref{algorithm_xi} and the definition of $\kappa_1$, we have $\left\|\Delta x_i(t)\right\|\leq 2\kappa_1\rho(t)$. Thus, by {Lemma \ref{lemma_varepsilon}}, we have
		\begin{align}\label{error_recursion_relation}
			\left\|e_i(t+1)\right\| \leq \varepsilon\left\|e_i(t)\right\|+2 \sqrt{N} \kappa_1 \rho(t),\quad \forall i \in \mathcal{N}
		\end{align}
		Based on the recursion relation \eqref{error_recursion_relation}, we have
		\begin{align*}
			||e_i(t)|| \leq\varepsilon^{t-1}||e_i(1)||+2 \sqrt{N} \kappa_1 \sum\nolimits_{k=0}^{t-2} \varepsilon^k \rho(t-k-1)
		\end{align*}
		This completes the proof. 
	\end{proof}
	
	Note that $0<\varepsilon<1$ and $\rho(t)$ is decreasing with $\lim\limits_{t\rightarrow\infty}\rho(t)=0$, we have $\lim\limits_{t\rightarrow\infty}\varepsilon^{t-1}=0$ and $\lim\limits_{t\rightarrow\infty}\sum_{k=0}^{t-2} \varepsilon^k \rho(t-k-1)$. Thus, $\lim\limits_{t\rightarrow\infty}||e_i(t)||=0$, and it implies the estimation will converge to the actual value.
	
	Then, a result on the bounds of the dual variables could be obtained.
	\begin{lemma}\label{lemma_dualvariable_bound}
		Under {Assumption} {\ref{assumption_domain_set}}, for $\forall i\in \mathcal{N}$ and $\forall t\in \mathcal{T}$,
		\begin{subequations}
			\begin{align}
				\sqrt{\rho(t)}\left\|\lambda_i(t)\right\| &\leq 3 \sqrt{N} \kappa_2\\
				\sqrt{\rho(t)}\left\|\mu_i(t)\right\| &\leq 3 \sqrt{N} \kappa_5
			\end{align}
		\end{subequations}
	\end{lemma}	
	Proof of {Lemma \ref{lemma_dualvariable_bound}} is similar to that of lemma 1 in \cite{Kaihong2021Environments}, which is omitted here due to the space limit.
	The next lemma provides an upper bound of the accumulated error between the v-GNE $x^*(t)$ and $x(t)$ obtained from Algorithm \ref{alg1}.
	\begin{lemma}\label{lemma_seekingGNE_bound}
		Under {Assumption} {\ref{assumption_domain_set}}, the accumulated error between $x^*(t)$ and $x(t)$ is upper bounded, i.e.,
		\begin{align}\label{eq_seekingGNE_bound}
			&\sum\nolimits_{t=1}^T\left\|x(t)-x^*(t)\right\|^2 \leq \nonumber\\
			&\quad\qquad\frac{2 \kappa_1}{\eta} \sum\nolimits_{t=1}^T \frac{1}{\rho^2(t)} \sum\nolimits_{i \in \mathcal{N}}\left\|x_i^*(t+1)-x_i^*(t)\right\|\nonumber\\
			&\quad\qquad+\frac{\pi_1}{2 \eta} \sum\nolimits_{t=1}^T \sqrt{\rho(t)}+\pi_2 \sum\nolimits_{t=1}^T \sum\nolimits_{i \in \mathcal{N}}\left\|e_i(t)\right\|\nonumber\\
			&\quad\qquad+\pi_3 \sum\nolimits_{t=1}^T \rho(t) \sum\nolimits_{i \in \mathcal{N}}\left\|e_i(t)\right\|^2\nonumber\\
			&\quad\qquad+\frac{1}{2 \eta} \sum\nolimits_{t=1}^T \frac{1}{\rho^2(t)}\left(\sum\nolimits_{i \in \mathcal{N}}\left\|x_i(t)-x_i^*(t)\right\|^2\right.\nonumber\\
			&\quad\qquad\left.-\sum\nolimits_{i \in\mathcal{N}}\left\|x_i(t+1)-x_i^*,(t+1)\right\|^2\right)
		\end{align}
		where 
		\begin{align*}
			&\pi_1= N\left(\delta_\lambda+\delta_\mu\right)^2+4N\left(\kappa_1+\kappa_3\right)\left(\delta_\lambda+\delta_\mu\right)+4 N \kappa_3^2\nonumber\\
			&\pi_2=2 \sqrt{N}(c+\theta)\left[2 \kappa_1+2 \kappa_3+\delta_\lambda+\delta_\mu\right]\nonumber\\
			&\pi_3=N c^2+\sqrt{N} c+\sqrt{N} c \theta^2+\theta^2
		\end{align*}
		with $\delta_\lambda=\kappa_4\big(3\sqrt{N}\kappa_2+\vartheta_\lambda\big)$, $\delta_\mu=\kappa_6\big(3\sqrt{N}\kappa_5+\vartheta_\mu\big)$, 
		$\vartheta_\lambda=\sup\limits_{i\in\mathcal{N},t\in\mathcal{T}}\left\|\lambda_i^*(t)\right\|$ and $\vartheta_\mu=\sup\limits_{i\in\mathcal{N},t\in\mathcal{T}}\left\|\mu_i^*(t)\right\|$.
	\end{lemma}	
	\begin{proof}Similar to (29) to (32) in\cite{Kaihong2021Environments}, by \eqref{algorithm_xi} and the definition of $\kappa_1$, we have the following two results.
		\begin{align}\label{traking_error}
			&\sum\limits_{i \in \mathcal{N}}\left\|x_i(t+1)-x_i^*(t+1)\right\|^2\leq \sum_{i \in \mathcal{N}}\left\|x_i(t+1)-x_i^*(t)\right\|^2 \nonumber\\
			&\qquad\qquad\qquad+4 \kappa_1 \sum\nolimits_{i \in \mathcal{N}}\left\|x_i^*(t+1)-x_i^*(t)\right\|\\
			\label{xi(t+1)-xi*(t)}
			&\sum\limits_{i \in \mathcal{N}}\left\|x_i(t+1)-x_i^*(t)\right\|^2
			\leq (1-\rho(t))\sum_{i \in \mathcal{N}}\left\|x_i(t)-x_i^*(t)\right\|^2\nonumber\\
			&\qquad\qquad+\rho(t)\sum\nolimits_{i \in \mathcal{N}}\left\|\mathcal{P}_{\chi_i}(\xi_i^1(t))-\mathcal{P}_{\chi_i}(\xi_i^2(t))\right\|^2
		\end{align}
		where 
		\begin{align*}
			&\xi_i^1(t)=x_i(t)-\rho(t)\left[\nabla_{x_i}J_{i,t}(\bar{x}^i(t))\right.\\
			&\ \left.+\rho(t)\left(A_i^T\lambda_i(t)+G_i^T\mu_i(t)\right)+c\sum\nolimits_{j \in \mathcal{N}_i}\left(x_i(t)-\bar{x}_i^j(t)\right)\right]\\
			&\xi_i^2(t)=x_i^*(t)-\rho(t)\left[\nabla_{x_i}J_{i,t}(x^*(t))\right.\nonumber\\
			&\qquad\qquad\qquad\qquad\qquad\left.+\rho(t)\left(A_i^T\lambda_i^*(t)+G_i^T\mu_i^*(t)\right)\right]
		\end{align*}
		For simplification, denote 
		\begin{align*}
			\phi_i^1&=x_i(t)-x_i^*(t),\
			\phi_i^2=\nabla_{x_i}J_{i,t}(\bar{x}^i(t))-\nabla_{x_i}J_{i,t}(x^*(t))\\ \phi_i^3&=A_i^T\lambda_i(t)-A_i^T\lambda_i^*(t),\  \phi_i^4=G_i^T\mu_i(t)-G_i^T\mu_i^*(t)\\ \phi_i^5&=c\sum\nolimits_{j \in \mathcal{N}_i}\big(x_i(t)-\bar{x}_i^j(t)\big), 
		\end{align*}
		and use $\rho$ to replace $\rho(t)$ in the remaining proof.
		By the non-expansive property of projection, we have
		\begin{align}\label{nonexpansive_property}
			&\left\|\mathcal{P}_{\chi_i}(\xi_i^1(t))-\mathcal{P}_{\chi_i}(\xi_i^2(t))\right\|^2
			\leq \left\|\xi_i^1(t)-\xi_i^2(t)\right\|^2\nonumber\\
			&=\left\|\phi_i^1\right\|^2+\rho^2\left\|\phi_i^2\right\|^2+\rho^4\left(\left\|\phi_i^3\right\|^2+\left\|\phi_i^4\right\|^2\right)+\rho^2\left\|\phi_i^5\right\|^2\nonumber\\
			&-2\rho\langle \phi_i^1,\phi_i^2\rangle-2\rho^2\langle \phi_i^1,\phi_i^3\rangle-2\rho^2\langle \phi_i^1,\phi_i^4\rangle-2\rho\langle \phi_i^1,\phi_i^5\rangle\nonumber\\
			&+2\rho^3\langle \phi_i^2,\phi_i^3\rangle+2\rho^3\langle \phi_i^2,\phi_i^4\rangle+2\rho^2\langle \phi_i^2,\phi_i^5\rangle\nonumber\\
			&+2\rho^4\langle \phi_i^3,\phi_i^4\rangle+2\rho^3\langle \phi_i^3,\phi_i^5\rangle+2\rho^3\langle \phi_i^4,\phi_i^5\rangle
		\end{align}
		
		Before continuing the proof, results on the bounds of the norm of $\phi_i^1, \phi_i^2, \cdots, \phi_i^5$ are given.
		\begin{align}
			\left\|\phi_i^1\right\|&\leq \left\| x_i(t)\right\|+\left\| x_i^*(t)\right\|\leq 2\kappa_1\\
			\left\|\phi_i^2\right\|&=\left\|\phi_i^{2,1}+\phi_i^{2,2}\right\|\leq \left\|\phi_i^{2,1}\right\|+2\kappa_3
		\end{align}
		where $\phi_i^{2,1}=\nabla_{x_i}J_{i,t}(\bar{x}^i(t))-\nabla_{x_i}J_{i,t}(x(t))$ and $\phi_i^{2,2}=\nabla_{x_i}J_{i,t}(x(t))-\nabla_{x_i}J_{i,t}(x^*(t))$.
		\begin{align}
			\sqrt{\rho}\left\|\phi_i^3\right\|&\leq\left\|A_i^T\left(\sqrt{\rho}\lambda_i(t)\right)\right\|+\sqrt{\rho}\left\|A_i^T\lambda_i^*(t)\right\|\nonumber\\
			&\leq 3\sqrt{N}\kappa_2\kappa_4+\kappa_4\vartheta_\lambda=\delta_\lambda
		\end{align}
		where the second inequality holds based on {Lemma \ref{lemma_dualvariable_bound}} and the fact that $\rho<1$. Similarly, we have
		\begin{align}
			\sqrt{\rho}\left\|\phi_i^4\right\|\leq\delta_\mu\label{phi_i^4}
		\end{align}
		By the definition of $e_i^j(t)$ and $e_i(t)$, we have
		\begin{align}
			\left\|\phi_i^5\right\|&=c\left\|\sum\nolimits_{j \in \mathcal{N}_i}e_i^j(t)\right\|\nonumber\\
			&\leq c\sum\nolimits_{j \in \mathcal{N}_i}\left\|e_i^j(t)\right\|\leq c\sqrt{N}\left\|e_i(t)\right\|
		\end{align}
		
		Thus, the sum of quadratic items in \eqref{nonexpansive_property} is upper bounded by
		\begin{align}\label{item22-55}
			&\rho^2\left\|\phi_i^2\right\|^2+\rho^4\left(\left\|\phi_i^3\right\|^2+\left\|\phi_i^4\right\|^2\right)+\rho^2\left\|\phi_i^5\right\|^2\nonumber\\
			&\leq \rho^2\left(\left\|\phi_i^{2,1}\right\|^2+4\kappa_3\left\|\phi_i^{2,1}\right\|+4\kappa_3^2\right)+\rho^3\left(\delta_\lambda^2+\delta_\mu^2\right)\nonumber\\
			&\quad+Nc^2\rho^2\left\|e_i(t)\right\|^2
		\end{align}
		Note that other items in \eqref{nonexpansive_property} are also upper bounded, i.e.,
		\begin{subequations}\label{item12-15}
			\begin{align}
				&-2\rho\langle\phi_i^1,\phi_i^2\rangle
				\leq 2\rho\left\|\phi_i^1\right\|\left\|\phi_i^{2,1}\right\|-2\rho\langle\phi_i^1,\phi_i^{2,2}\rangle\nonumber\\
				&\qquad\qquad\leq 4\rho\kappa_1\left\|\phi_i^{2,1}\right\|-2\rho\langle\phi_i^1,\phi_i^{2,2}\rangle\label{item12}\\
				&-2\rho^2\langle\phi_i^1,\phi_i^3\rangle-2\rho^2\langle\phi_i^1,\phi_i^4\rangle\nonumber\\
				&\qquad\qquad\leq 2\rho^{\frac{3}{2}}\left\|\phi_i^1\right\|\left(\sqrt{\rho}\left\|\phi_i^3\right\|+\sqrt{\rho}\left\|\phi_i^4\right\|\right)\nonumber\\
				&\qquad\qquad\leq 4\kappa_1\rho^{\frac{3}{2}}\left(\delta_\lambda+\delta_\mu\right)\label{item13-14}\\
				&-2\rho\langle\phi_i^1,\phi_i^5\rangle\leq2\rho\left\|\phi_i^1\right\|\left\|\phi_i^5\right\|\leq     4c\sqrt{N}\kappa_1\rho\left\|e_i(t)\right\|\label{item15}
			\end{align}
		\end{subequations}
		Similarly, we have
		\begin{subequations}\label{item23-45}
			\begin{align}
				2\rho^3\langle\phi_i^2,\phi_i^3\rangle&\leq 2\delta_\lambda\rho^2\left(\left\|\phi_i^{2,1}\right\|+2\kappa_3\right)\label{item23}\\
				2\rho^3\langle\phi_i^2,\phi_i^4\rangle&\leq 2\delta_\mu\rho^2\left(\left\|\phi_i^{2,1}\right\|+2\kappa_3\right)\label{item24}\\
				2\rho^2\langle\phi_i^2,\phi_i^5\rangle&\leq 2c\sqrt{N}\rho^2\left\|e_i(t)\right\|\left(\left\|\phi_i^{2,1}\right\|+2\kappa_3\right)\label{item25}\\
				2\rho^4\langle\phi_i^3,\phi_i^4\rangle&\leq 2\delta_\lambda\delta_\mu\rho^3\label{item34}\\
				2\rho^3\langle\phi_i^3,\phi_i^5\rangle&\leq 2c\sqrt{N}\delta_\lambda\rho^2\left\|e_i(t)\right\|\label{item35}\\
				2\rho^3\langle\phi_i^4,\phi_i^5\rangle&\leq 2c\sqrt{N}\delta_\mu\rho^2\left\|e_i(t)\right\|\label{item45}
			\end{align}
		\end{subequations}
		
		Substituting \eqref{item22-55}, \eqref{item12-15}, \eqref{item23-45} into \eqref{xi(t+1)-xi*(t)}, we have
		\begin{align}\label{mid_item}
			&\sum\nolimits_{i \in \mathcal{N}}\left\|x_i(t+1)-x_i^*(t)\right\|^2\nonumber\\
			\leq&\sum\nolimits_{i \in \mathcal{N}}\left\|x_i(t)-x_i^*(t)\right\|^2+N\left(\delta_\lambda+\delta_\mu\right)^2\rho^4\nonumber\\
			&+4N\kappa_3\left(\delta_\lambda+\delta_\mu+\kappa_3\right)\rho^3+4N\kappa_1\left(\delta_\lambda+\delta_\mu\right)\rho^{\frac{5}{2}}\nonumber\\
			&+2\sqrt{N}c\rho^2\left[2\kappa_1+\left(2\kappa_3+\delta_\lambda+\delta_\mu\right)\rho\right]\sum\nolimits_{i \in \mathcal{N}}\left\|e_i(t)\right\|\nonumber\\
			&+Nc^2\rho^3\sum\nolimits_{i \in \mathcal{N}}\left\|e_i(t)\right\|^2\nonumber\\
			&-2\rho^2\langle\phi_i^1,\phi_i^{2,2}\rangle+\rho^3\sum\nolimits_{i \in \mathcal{N}}\left\|\phi_i^{2,1}\right\|^2\nonumber\\
			&+2\rho^2
			\left[2\kappa_1+\left(2\kappa_3+\delta_\lambda+\delta_\mu\right)\rho\right]\sum\nolimits_{i \in \mathcal{N}}\left\|\phi_i^{2,1}\right\|\nonumber\\
			&+2\sqrt{N}c\rho^3\sum\nolimits_{i \in \mathcal{N}}\left\|e_i(t)\right\|\left\|\phi_i^{2,1}\right\|
		\end{align}
		
		Then, by {Lemma \ref{lemma_monotonicity_continuity}}, results on the upper bounds of the last four items in \eqref{mid_item} are given as follows.
		
		Due to the Lipschitz continuity of $F_t(x)$, we have
		\begin{align}\label{item_sum_1}
			&\sum\nolimits_{i \in \mathcal{N}}\left\|\phi_i^{2,1}\right\|\leq\theta\sum\nolimits_{i \in \mathcal{N}}\left\|\bar{x}^i(t)-x(t)\right\|\nonumber\\
			\leq& \theta\sqrt{\bigg(\sum_{i \in \mathcal{N}}\left\|\bar{x}^i(t)-x(t)\right\|\bigg)^2}\leq\theta\sqrt{N\sum_{i \in \mathcal{N}}\left\|\bar{x}^i(t)-x(t)\right\|^2}\nonumber\\
			=&\theta\sqrt{N\sum\nolimits_{i \in \mathcal{N}}\left\|e_i(t)\right\|^2}\leq \theta\sqrt{N}\sum\nolimits_{i \in \mathcal{N}}\left\|e_i(t)\right\|
		\end{align}
		and
		\begin{align}\label{item_sum_2}
			\sum\nolimits_{i \in \mathcal{N}}\left\|\phi_i^{2,1}\right\|^2\leq \theta^2\sum\nolimits_{i \in \mathcal{N}}\left\|e_i(t)\right\|^2
		\end{align}
		Based on \eqref{item_sum_2}, we have
		\begin{align}\label{item_sum_3}
			\sum\nolimits_{i \in \mathcal{N}}\left\|e_i(t)\right\|\left\|\phi_i^{2,1}\right\|\leq&\frac{1}{2}\sum\nolimits_{i \in \mathcal{N}}\left(\left\|e_i(t)\right\|^2+\left\|\phi_i^{2,1}\right\|^2\right)\nonumber\\
			\leq&\frac{1}{2}\left(1+\theta^2\right)\sum\nolimits_{i \in \mathcal{N}}\left\|e_i(t)\right\|^2
		\end{align}
		Moreover, since $F_t(x)$ is $\eta-$strongly monotone, we have
		\begin{align}\label{item_sum_4}
			-\sum\nolimits_{i \in \mathcal{N}}\langle \phi_i^1,\phi_i^{2,2}\rangle
			&=-\langle x(t)-x^*(t),F_t(x(t))-F_t(x^*(t))\rangle\nonumber\\
			&\leq-\eta\left\|x(t)-x^*(t)\right\|^2
		\end{align}
		
		Substituting $\eqref{item_sum_1} \sim \eqref{item_sum_4}$ into \eqref{mid_item}, we have
		\begin{align}\label{pre_of_lemma6}
			&\sum\nolimits_{i \in \mathcal{N}}\left\|x_i(t+1)-x_i^*(t)\right\|^2
			\leq\sum\nolimits_{i \in \mathcal{N}}\left\|x_i(t)-x_i^*(t)\right\|^2\nonumber\\
			&-2\eta\rho^2\left\|x(t)-x^*(t)\right\|^2+\pi_1\rho^{\frac{5}{2}}+\pi_2\rho^2\sum\nolimits_{i \in \mathcal{N}}\left\|e_i(t)\right\|\nonumber\\
			&+\pi_3\rho^3\sum\nolimits_{i \in \mathcal{N}}\left\|e_i(t)\right\|^2
		\end{align}
		% 		where $\pi_1$, $\pi_2$ and $\pi_3$ are defined in \eqref{pai_definition}.
		Substituting \eqref{pre_of_lemma6} into \eqref{traking_error} yields \eqref{eq_seekingGNE_bound}.
	\end{proof}
	
	Finally, we can prove {Theorem \ref{Theorem_1}}.
	
	\begin{proof}
		Following \cite{Kaihong2021Environments}, by {Lemma \ref{lemma_error_bound}}, we have
		\begin{subequations}\label{corollary_of_lemma4}
			\begin{align}
				\sum\nolimits_{t=1}^T \sum\nolimits_{i \in \mathcal{N}}\left\|e_i(t)\right\| &\leq \mathcal{O}\left(\sum\nolimits_{t=1}^T \rho(t)\right) \\
				\sum_{t=1}^T \rho(t) \sum_{i \in \mathcal{N}}\left\|e_i(t)\right\|^2 &\leq \mathcal{O}\left(\sum_{t=1}^T \rho(t)\right)
			\end{align}
		\end{subequations}
		Using the fact that $\rho(t+1) \leq \rho(t)$ and $\left\|x_i(t)-x_i^*(t)\right\|^2 \leq\left(\left\|x_i(t)\right\|+\left\|x_i^*(t)\right\|\right)^2 \leq 4 \kappa_1^2$, we have
		\begin{align}\label{bound_of_the_last_item_of_lemma6}
			&\sum\nolimits_{t=1}^T \frac{1}{\left(\rho(t)\right)^2}\left(\sum\nolimits_{i \in \mathcal{N}}\left\|x_i(t)-x_i^*(t)\right\|^2\right.\nonumber\\
			&\qquad\quad\left.-\sum\nolimits_{i \in \mathcal{N}}\left\|x_i(t+1)-x_i^*(t+1)\right\|^2\right)\nonumber\\
			=&\frac{1}{\rho^2(1)} \sum\nolimits_{i \in \mathcal{N}}\left\|x_i(1)-x_i^*(1)\right\|^2\nonumber\\
			&-\frac{1}{\rho^2(T)} \sum\nolimits_{i \in \mathcal{N}}\left\|x_i(T+1)-x_i^*(T+1)\right\|^2\nonumber\\
			&+\sum\nolimits_{t=2}^T\left[\frac{1}{\rho^2(t)}-\frac{1}{\rho^2(t-1)}\right] \sum\nolimits_{i \in \mathcal{N}}\left\|x_i(t)-x_i^*(t)\right\|^2\nonumber\\
			%\leq& \frac{1}{\rho^2{(1)}} \sum\nolimits_{i \in \mathcal{N}}\left\|x_i(1)-x_i^*(1)\right\|^2\nonumber\\
			%&-\frac{1}{\rho^2{(T)}} \sum\nolimits_{i \in \mathcal{N}}\left\|x_i(T+1)-x_i^*(T+1)\right\|^2\nonumber\\
			%&+\sum\nolimits_{t=2}^T 4 N \kappa_1^2\left[\frac{1}{\rho^2{(t)}}-\frac{1}{\rho^2(t-1)}\right]\nonumber\\
			%=&\frac{1}{\rho^2(T)}\left(4 N \kappa_1^2-\sum\nolimits_{i \in \mathcal{N}}\left\|x_i(T+1)-x_i^*(T+1)\right\|^2\right)\nonumber\\
			%&-\frac{1}{\rho^2(1)}\left(4 N \kappa_1^2-\sum\nolimits_{i \in \mathcal{N}}\left\|x_i(1)-x_i^*(1)\right\|^2\right)\nonumber\\
			\leq& \frac{4 N \kappa_1^2}{\rho^2(T)}
		\end{align}
		Moreover, by the definition of $\varPhi_T$ and the monotonic descending of $\rho(t)$, we have
		\begin{align}\label{bound_of_the_first_item_of_lemma6}
			\sum\nolimits_{t=1}^T &\frac{1}{\rho^2(t)} \sum\nolimits_{i \in \mathcal{N}}\left\|x_i^*(t+1)-x_i^*(t)\right\|\nonumber\\
			\leq& \frac{1}{\rho^2(T)} \sum\nolimits_{t=1}^T \sum\nolimits_{i \in \mathcal{N}}\left\|x_i^*(t+1)-x_i^*(t)\right\|\nonumber\\
			\leq& \frac{1}{\rho^2(T)} \sum\nolimits_{t=1}^T \sqrt{N}\left\|x^*(t+1)-x^*(t)\right\|\nonumber\\
			=&\frac{\sqrt{N} \varPhi_T}{\rho^2(T)}
		\end{align}
		Substituting $\eqref{corollary_of_lemma4}$,  $\eqref{bound_of_the_last_item_of_lemma6}$ and $\eqref{bound_of_the_first_item_of_lemma6}$ into $\eqref{eq_seekingGNE_bound}$ yields
		\begin{align}\label{bound_of_sequence_error}
			&\sum\nolimits_{t=1}^T\left\|x(t)-x^*(t)\right\|^2 
			\leq \frac{2 \sqrt{N} \kappa_1}{\eta\rho^2(T)}\left(\varPhi_T+\sqrt{N} \kappa_1\right)\nonumber\\
			&+\frac{\pi_1}{2 \eta} \sum\nolimits_{t=1}^T \sqrt{\rho(t)}+\pi_2 \sum\nolimits_{t=1}^T \sum\nolimits_{i \in \mathcal{N}}\left\|e_i(t)\right\| \nonumber\\
			&+\pi_3 \sum\nolimits_{t=1}^T \rho(t) \sum\nolimits_{i \in \mathcal{N}}\left\|e_i(t)\right\|^2 \nonumber\\
			&\leq  \mathcal{O}\left(\frac{\varPhi_T+1}{\rho^2(T)}+\sum\nolimits_{t=1}^T \sqrt{\rho(t)}\right)
		\end{align}
		Therefore, by the definition of $\kappa_3$, we have
		\begin{align}\label{bound_of_regret}
			R_i(T) &=\sum\nolimits_{t=1}^T\left(J_{i,t}\left(x_i(t), x_{-i}^*(t)\right)-J_{i,t}\left(x^*(t)\right)\right) \nonumber\\
			& \leq \kappa_3 \sum\nolimits_{t=1}^T\left\|x_i(t)-x^*(t)\right\| \nonumber\\
			& \leq \kappa_3 \sum\nolimits_{t=1}^T\left\|x(t)-x^*(t)\right\| \nonumber\\
			& \leq \kappa_3 \sqrt{T \sum\nolimits_{t=1}^T\left\|x(t)-x^*(t)\right\|^2} \nonumber\\
			& \leq \mathcal{O}\left(\sqrt{T\left(\frac{\varPhi_T+1}{\rho^2(T)}+\sum_{t=1}^T\sqrt{\rho(t)}\right)}\right)
		\end{align}
		This completes the proof. 
	\end{proof}

	\bibliographystyle{IEEEtran}
	\bibliography{mybib}

\begin{thebibliography}{10}
\providecommand{\url}[1]{#1}
\csname url@samestyle\endcsname
\providecommand{\newblock}{\relax}
\providecommand{\bibinfo}[2]{#2}
\providecommand{\BIBentrySTDinterwordspacing}{\spaceskip=0pt\relax}
\providecommand{\BIBentryALTinterwordstretchfactor}{4}
\providecommand{\BIBentryALTinterwordspacing}{\spaceskip=\fontdimen2\font plus
\BIBentryALTinterwordstretchfactor\fontdimen3\font minus
  \fontdimen4\font\relax}
\providecommand{\BIBforeignlanguage}[2]{{%
\expandafter\ifx\csname l@#1\endcsname\relax
\typeout{** WARNING: IEEEtran.bst: No hyphenation pattern has been}%
\typeout{** loaded for the language `#1'. Using the pattern for}%
\typeout{** the default language instead.}%
\else
\language=\csname l@#1\endcsname
\fi
#2}}
\providecommand{\BIBdecl}{\relax}
\BIBdecl

\bibitem{chen2019energy}
Y.~Chen, S.~Mei, F.~Zhou, S.~H. Low, W.~Wei, and F.~Liu, ``An energy sharing
  game with generalized demand bidding: Model and properties,'' \emph{IEEE
  Trans. Smart Grid}, vol.~11, no.~3, pp. 2055--2066, 2019.

\bibitem{Wang2021Distributed}
Z.~Wang, F.~Liu, Z.~Ma, Y.~Chen, M.~Jia, W.~Wei, and Q.~Wu, ``Distributed
  generalized nash equilibrium seeking for energy sharing games in prosumers,''
  \emph{IEEE Trans. Power Syst.}, vol.~36, no.~5, pp. 3973--3986, 2021.

\bibitem{Morstyn2018peer-to-peer}
T.~Morstyn, N.~Farrell, S.~J. Darby, and M.~D. Mcculloch, ``Using peer-to-peer
  energy-trading platforms to incentivize prosumers to form federated power
  plants,'' \emph{Nature Energy}, vol.~3, no.~2, p. 94–101, 2018.

\bibitem{Cui2020FairPeer-to-Peer}
S.~Cui, Y.-W. Wang, Y.~Shi, and J.-W. Xiao, ``A new and fair peer-to-peer
  energy sharing framework for energy buildings,'' \emph{IEEE Trans. Smart
  Grid}, vol.~11, no.~5, pp. 3817--3826, 2020.

\bibitem{Nespoli2018rational}
L.~Nespoli, M.~Salani, and V.~Medici, ``A rational decentralized generalized
  nash equilibrium seeking for energy markets,'' in \emph{2018 International
  Conference on Smart Energy Systems and Technologies (SEST)}, 2018, pp. 1--6.

\bibitem{Chen2021SameFlexibility}
Y.~Chen, W.~Wei, H.~Wang, Q.~Zhou, and J.~P.~S. Catalão, ``An energy sharing
  mechanism achieving the same flexibility as centralized dispatch,''
  \emph{IEEE Trans. Smart Grid}, vol.~12, no.~4, pp. 3379--3389, 2021.

\bibitem{Belgioioso2022Operationally}
G.~Belgioioso, W.~Ananduta, S.~Grammatico, and C.~Ocampo-Martinez,
  ``Operationally-safe peer-to-peer energy trading in distribution grids: A
  game-theoretic market-clearing mechanism,'' \emph{IEEE Trans. Smart Grid},
  vol.~13, no.~4, pp. 2897--2907, 2022.

\bibitem{Hans2019Hierarchical}
C.~A. Hans, P.~Braun, J.~Raisch, L.~Grüne, and C.~Reincke-Collon,
  ``Hierarchical distributed model predictive control of interconnected
  microgrids,'' \emph{IEEE Trans. Sustainable Energy}, vol.~10, no.~1, pp.
  407--416, 2019.

\bibitem{Atzeni2013Demand-Side}
I.~Atzeni, L.~G. Ordóñez, G.~Scutari, D.~P. Palomar, and J.~R. Fonollosa,
  ``Demand-side management via distributed energy generation and storage
  optimization,'' \emph{IEEE Trans. Smart Grid}, vol.~4, no.~2, pp. 866--876,
  2013.

\bibitem{Facchinei2010Generalized}
F.~Facchinei and C.~Kanzow, ``Generalized {N}ash equilibrium problems,''
  \emph{Annals of Oper. Res.}, vol. 175, no.~1, pp. 177--211, 2010.

\bibitem{Kulkarni2012variational}
A.~Kulkarni and U.~Shanbhag, ``On the variational equilibrium as a refinement
  of the generalized nash equilibrium,'' \emph{Automatica}, vol.~48, no.~1, pp.
  45--55, 2012.

\bibitem{Facc2010Nashequilibria}
F.~Facchinei and J.-S. Pang, ``Nash equilibria: The variational approach,''
  \emph{Convex Optimization in Signal Processing and Communications}, pp.
  443--449, 2010.

\bibitem{wang2019distributed2}
Z.~Wang, F.~Liu, S.~H. Low, C.~Zhao, and S.~Mei, ``Distributed frequency
  control with operational constraints, part {II}: Network power balance,''
  \emph{IEEE Trans. Smart Grid}, vol.~10, no.~1, pp. 53--64, Jan 2019.

\bibitem{wang2019distributed}
Z.~Wang, F.~Liu, J.~Z. Pang, S.~H. Low, and S.~Mei, ``Distributed optimal
  frequency control considering a nonlinear network-preserving model,''
  \emph{IEEE Trans. Power Syst.}, vol.~34, no.~1, pp. 76--86, 2019.

\bibitem{wang2019distributed_variation}
Z.~Wang, S.~Mei, F.~Liu, S.~H. Low, and P.~Yang, ``Distributed load-side
  control: Coping with variation of renewable generations,'' \emph{Automatica},
  vol. 109, p. 108556, 2019.

\bibitem{pavel2019distributed}
L.~{Pavel}, ``Distributed gne seeking under partial-decision information over
  networks via a doubly-augmented operator splitting approach,'' \emph{IEEE
  Trans. Autom. Control}, vol.~65, no.~4, pp. 1584--1597, 2020.

\bibitem{Kaihong2021Environments}
K.~Lu, G.~Li, and L.~Wang, ``Online distributed algorithms for seeking
  generalized nash equilibria in dynamic environments,'' \emph{IEEE Trans.
  Autom. Control}, vol.~66, no.~5, pp. 2289--2296, 2021.

\bibitem{2021arXiv210506200M}
M.~{Meng}, X.~{Li}, Y.~{Hong}, J.~{Chen}, and L.~{Wang}, ``{Decentralized
  Online Learning for Noncooperative Games in Dynamic Environments},''
  \emph{arXiv e-prints}, p. arXiv:2105.06200, May 2021.

\bibitem{Shahrampour2018MirrorDescent}
S.~Shahrampour and A.~Jadbabaie, ``Distributed online optimization in dynamic
  environments using mirror descent,'' \emph{IEEE Trans. Autom. Control},
  vol.~63, no.~3, pp. 714--725, 2018.

\bibitem{Hall2015Online}
E.~C. Hall and R.~M. Willett, ``Online convex optimization in dynamic
  environments,'' \emph{IEEE J. Sel. Top. Signal Process.}, vol.~9, no.~4, pp.
  647--662, 2015.

\bibitem{Kaihong2020Pseudoconvex}
K.~Lu, G.~Jing, and L.~Wang, ``Online distributed optimization with strongly
  pseudoconvex-sum cost functions,'' \emph{IEEE Trans. Autom. Control},
  vol.~65, no.~1, pp. 426--433, 2021.

\bibitem{website_PV}
NREL. (2018) Solrmap utah geological survey.
  \url{https://midcdmz.nrel.gov/usep_cedar/}.

\bibitem{website_load}
UCI. (2012) Individual household electric power consumption data set.
  \url{https://archive.ics.uci.edu/ml/datasets/individual+household+electric+power+consumption}.

\bibitem{Mesbahi2010Graph}
M.~Mesbahi and M.~Egerstedt, \emph{Graph Theoretic Methods for Multiagent
  Networks}.\hskip 1em plus 0.5em minus 0.4em\relax Princeton: Princeton
  University Press, 2010.

\end{thebibliography}
\end{document}